\documentclass[12pt,twoside,doublespace]{article}

\usepackage{latexsym,amssymb,amsfonts}
\usepackage{graphicx}

\def\smskip{\par\vskip 5 pt}
\def\QED{\hfill $\Box$\smskip}
\newtheorem{theorem}{Theorem}
\newtheorem{lemma}{Lemma}
\newtheorem{proposition}{Proposition}
\newtheorem{corollary}{Corollary}
\newtheorem{definition}{Definition}

\begin{document}

\begin{center}

\vspace{35pt}

{\Large \bf Primal-Dual Method for Optimization Problems}

\vspace{5pt}

{\Large \bf  with Changing Constraints}

\vspace{5pt}

\vspace{35pt}

{\sc Igor V.~Konnov\footnote{\normalsize E-mail: konn-igor@ya.ru}}

\vspace{35pt}

{\em  Department of System Analysis
and Information Technologies, \\ Kazan Federal University, ul.
Kremlevskaya, 18, Kazan 420008, Russia.}

\end{center}

\vspace{25pt}

\hrule

\vspace{10pt}

{\bf Abstract:} We propose a modified primal-dual method
for general convex optimization problems with changing constraints.
We obtain properties of Lagrangian saddle points for these problems
which enable us to establish convergence of the proposed method.
We describe specializations of the proposed approach
to  multi-agent optimization problems under changing communication topology
and to feasibility problems.

{\bf Key words:} Convex optimization, changing constraints, primal-dual method,
constrained multi-agent optimization, feasibility problem.

\vspace{10pt}

\hrule

\vspace{10pt}

{\em MS Classification:} {65K05, 90C06, 90C25, 68M14, 68W15, 93A14}


\section{Introduction}\label{sc:1}

It is well known that the general optimization problem consists in
finding the minimal value of some goal function $\tilde f$ on
a feasible set $\tilde D$. For brevity, we write
this problem as
$$
 \min \limits _{v \in \tilde D} \to \tilde f (v).
$$
In many cases, only some approximations are known instead of the
exact values of the goal function and the feasible set. This
situation is caused by various circumstances. On the one hand, this
is due to inevitable calculation errors of values of cost and
constraint functions. On the other hand, this is due to
incompleteness of information about these functions since their
parameters may be specialized during the computational process. Such
problems are called non stationary; see e.g. \cite{EM79} and  \cite[Chapter VI, \S
3]{Pol83}. Besides, some perturbations can be inserted for attaining better properties in comparison with the initial one as in
various regularization methods; see e.g. \cite{Vas81}.
In these problems,  only some sequences of approximations $\{ \tilde D_{k} \}$
and $\{ \tilde f_{k} \}$ are known, which however must converge in some sense to
the exact values of $\tilde D$ and  $\tilde f$.
The case where the convergence is not obligatory seems more difficult, but
it also appears in many applied problems. For instance, large-scale models may contain
superfluous constraints and variables together with the necessary ones, but only some of them can be utilized
at a given iterate.  Various decentralized multi-agent optimization problems can serve as examples of
such systems; see e.g. \cite{KPK09,LOF11,PYY13} and the references therein.

In this paper we investigate just general  convex optimization problems with changing constraints.
First we obtain properties of Lagrangian saddle points for these problems.
They enable us to propose a modification of the
primal-dual method from \cite{Ant88} for finding their solutions.
We establish different convergence properties of the proposed method
under rather weak assumptions. We describe specializations of the proposed approach
to  multi-agent optimization problems under changing communication topology
and to feasibility problems.


\section{The general problem with changing constraints and its properties}\label{sc:2}

Let us consider first a general optimization problem of the form
\begin{equation} \label{eq:2.1}
 \min \limits _{x \in D} \to f (x)
\end{equation}
for some  function $f : \mathbb{E} \to \mathbb{R}$ and set $D \subseteq \mathbb{E}$
in a finite-dimensional space $\mathbb{E}$.
The set of its solutions is denoted by $D^{*}$, and the optimal
function value by $f^{*}$, i.e.
$$
f^{*} = \inf \limits_{ x \in D} f(x).
$$
It will be suitable for us to specialize this problem as follows.
For each $x \in \mathbb{E}$, let $x=(x_{i})_{i=1, \ldots, m}$, i.e. $x^{\top}=(x^{\top}_{1}, \dots,x^{\top}_{m})$,
where $x_{i}= (x_{i1}, \dots,x_{in})^{\top}$ for $i=1, \dots, m$,
hence $\mathbb{E}=\mathbb{R}^{mn}$. This means that each vector $x$ is divided into $m$ subvectors
$x_{i} \in \mathbb{R}^{n}$. In case $n=1$ we obtain the custom coordinates of $x$.
Next, we suppose that
\begin{equation} \label{eq:2.2}
 D = \left\{x \in X \ | \ A x=b  \right\},
\end{equation}
where $X$ is a subset of $\mathbb{R}^{mn}$, the matrix $A$ has $ln$ rows and $mn$ columns, so that
$b =(b_{i})_{i=1, \ldots, l}$, $b_{i} \in \mathbb{R}^{n}$ for $i=1, \dots, m$, and $b \in \mathbb{R}^{ln}$.

In what follows, we will use the following basic assumptions.
\begin{enumerate}
\item[{\bf (A1)}] The set $ D^{*}$ is nonempty, $X$ is a  convex and closed set in
$ \mathbb{R}^{mn}$.
\item[{\bf (A2)}] $f : \mathbb{R}^{mn} \to \mathbb{R}$ is a convex function.
\end{enumerate}

For brevity, we set $M = \{1, \dots, m\}$ and $L = \{1, \dots, l\}$.
It is clear that the matrix $A$ is represented as follows:
$$
A=
\left( {
\begin{array}{c}
A_{1} \\
A_{2} \\
\dots \\
A_{l}
\end{array}
} \right),
$$
where $A_{i}$ is the corresponding $n \times mn$ sub-matrix of $A$ for $i \in L$.
We will write this briefly
$$
A=\left( \{A^{\top}_{i}\} _{i \in L}\right)^{\top}.
$$
Similarly, we can determine some other submatrices
$$
A_{I}=\left( \{A^{\top}_{i}\} _{i \in I}\right)^{\top}
$$
for any $I \subseteq L$, hence $A=A_{L}$. Setting
\begin{equation} \label{eq:2.3}
 F_{I} = \left\{x \in \mathbb{R}^{mn} \ | \ A_{I} x=b_{I}  \right\} \ \mbox{and} \
 D_{I} = \left\{x \in X \ | \ A_{I} x=b_{I}  \right\}= X \bigcap  F_{I},
\end{equation}
where $b_{I}=(b_{i})_{i \in I}$, we obtain a family of
optimization problems
\begin{equation} \label{eq:2.4}
 \min \limits _{x \in D_{I}} \to f (x).
\end{equation}
As above, we denote the solution set of problem (\ref{eq:2.3})--(\ref{eq:2.4}) by $D^{*}_{I}$, and the optimal
function value by $f^{*}_{I}$, so that $D^{*}_{L}=D^{*}$ and $f^{*}_{L}=f^{*}$.
Clearly, if  $I \subset J$, then $f^{*}_{I} \leq f^{*}_{J}$.
We intend to establish some properties related to superfluous constraints.
We will denote by $F^{*}$ the solution set of the optimization problem
$$
 \min \limits _{x \in X} \to f (x),
$$
and its optimal
function value by $f^{**}$.

\begin{lemma} \label{lm:2.1} Suppose  the set $F^{*} \bigcap  F_{I} $ is nonempty for some
$I \subseteq L$. Then $f^{**} = f^{*}_{I}$  and $F^{*} \bigcap  F_{I} = D^{*}_{I}$.
\end{lemma}
{\bf Proof.} If $x^{*} \in F^{*} \bigcap  F_{I}$, then clearly $x^{*} \in D^{*}_{I}$, hence
$f^{**} = f^{*}_{I}$. It follows that $F^{*} \bigcap  F_{I} = D^{*}_{I}$.
\QED


\begin{definition} \label{def:2.1}
We say that $I \subset J $ is a basic index set with respect  to $J$  if
$$
A_{I} x=b_{I} \ \Longrightarrow \ A_{J} x=b_{J}.
$$
We say that $I \subset L$ is a basic index set if it is a basic index set
with respect to $L$.
\end{definition}

From the definitions we obtain immediately the simple but useful properties.

\begin{lemma} \label{lm:2.2} \hfill
\begin{description}
  \item{(i)} If $I \subset J $ is a basic index set with respect to $J$, then
$f^{*}_{I} = f^{*}_{J}$, $D_{I} = D_{J}$, and $D^{*}_{I} = D^{*}_{J}$.
  \item{(ii)} If $I $ is a basic index set, then $f^{*}_{I} = f^{*}$, $D_{I} = D$,  and $D^{*}_{I} = D^{*}$.
\end{description}
\end{lemma}

For each problem (\ref{eq:2.3})--(\ref{eq:2.4}) associated with an index set $I \subseteq L$
we can define its Lagrange function
$$
 \mathcal{L}_{I}(x,y) =  f(x)+  \langle y_{I}, A_{I} x-b_{I} \rangle
$$
and the corresponding saddle point problem. It appears more suitable to utilize the general
Lagrange function
$$
 \mathcal{L}(x,y) =  f(x)+  \langle y, A x-b \rangle,
$$
with the modified dual feasible set. Namely, we say that $w^{*}=(x^{*},y^{*}) \in X\times
Y_{I}$ is a saddle point for problem (\ref{eq:2.3})--(\ref{eq:2.4}) if
\begin{equation} \label{eq:2.5}
 \forall y \in Y_{I}, \quad  \mathcal{L}(x^*,y)\leq \mathcal{L}(x^{*},y^{*})\leq \mathcal{L}(x,y^{*}) \quad \forall x\in X,
\end{equation}
where
$$
Y_{I}=\left\{y=(y_{i})_{i \in L}\in \mathbb{R}^{ln} \ | \ y_{i}=\mathbf{0} \in \mathbb{R}^{n} \ \mbox{for} \  i \notin I  \right\}.
$$
We denote by $ W^{*}_{I}=D^{*}_{I}\times Y^{*}_{I}$ the set of saddle points in (\ref{eq:2.5}) since
$D^{*}_{I}$ is precisely the solution set of problem (\ref{eq:2.3})--(\ref{eq:2.4}), whereas
$ Y^{*}_{I}$ is the set of its Lagrange multipliers. Since $D^{*}_{L}=D^{*}$, we also set
$ Y^{*}=Y^{*}_{L}$, i.e. $W^{*}= D^{*}\times Y^{*}$ is the set of saddle points for the
initial problem (\ref{eq:2.1})--(\ref{eq:2.2}).
Observe that (\ref{eq:2.5}) is rewritten equivalently as follows:
\begin{equation} \label{eq:2.6}
 A_{I} x^*=b_{I}, \quad  \mathcal{L}(x^{*},y^{*})\leq \mathcal{L}(x,y^{*}) \quad \forall x\in X.
\end{equation}
Besides, if we take $I =\varnothing $, then $ Y_{I}=\{\mathbf{0} \}$, hence we can write
$ D^{*}_{I}= F^{*}$ and $ Y^{*}_{I}=\{\mathbf{0} \}$.


\begin{proposition} \label{pro:2.1}
Suppose that assumptions (A1)--(A2) are fulfilled.
If $I \subset J $ is a basic index set with respect to $J$, then
$D^{*}_{I} = D^{*}_{J}$ and $Y^{*}_{I} \subseteq Y^{*}_{J}$.
\end{proposition}
{\bf Proof.} The first equality follows from Lemma \ref{lm:2.2} (i).
If $(x^{*},y^{*}) \in D^{*}_{I}\times Y^{*}_{I}$, then (\ref{eq:2.6}) holds,
which now implies (\ref{eq:2.6}) with $I=J$.
Hence $y^{*} \in Y^{*}_{J}$.
\QED


\begin{corollary} \label{cor:2.1}
Suppose that assumptions (A1)--(A2) are fulfilled.
If $I $ is a basic index set, then $D^{*}_{I} = D^{*}$ and $Y^{*}_{I} \subseteq Y^{*}$.
\end{corollary}

We can establish similar relations for dual variables in case $F^{*} \bigcap  F_{I} \neq \varnothing $.


\begin{proposition} \label{pro:2.2}
Suppose that assumptions (A1)--(A2) are fulfilled,  the set $F^{*} \bigcap  F_{I} $ is nonempty for some
$I \subseteq L$. Then $F^{*} \bigcap  F_{I} = D^{*}_{I}$  and $\mathbf{0} \in  Y^{*}_{I}$.
\end{proposition}
{\bf Proof.} The first equality follows from Lemma \ref{lm:2.1}.
Take any $x^{*} \in F^{*} \bigcap  F_{I}$, then $x^{*} \in D^{*}_{I}$
and (\ref{eq:2.6})  holds with $y^{*}=\mathbf{0}$.
Therefore,  $\mathbf{0} \in  Y^{*}_{I}$.
\QED


\section{Primal-dual method for the family of saddle point problems}\label{sc:3}

We intend to find saddle points in (\ref{eq:2.5}) by a modification of the
primal-dual method that was proposed in \cite{Ant88}.
First we note that the set of saddle points for the
initial problem (\ref{eq:2.1})--(\ref{eq:2.2}) is nonempty under the assumptions in
(A1)--(A2); see e.g. \cite[Corollary 28.2.2]{Roc70}. Therefore, this is the case
for each saddle point problem in (\ref{eq:2.5}) associated with a
basic index set $I $. Denote by $\pi _{U} (u)$ the
projection of $u$ onto $U$. Also, for simplicity we will write $Y_{(k)}=Y_{I_{k}}$,
$Y^{*}_{(k)}=Y^{*}_{I_{k}}$, etc. Then the method is described as follows.

\medskip
\noindent {\bf Method (PDM).}
{\em Step 0:} Choose an index set $I_{0} \subseteq L$, a point
$w^{0}=(x^{0},y^{0}) \in X\times Y_{(0)}$. Set $k=1$.

{\em Step 1:} Choose an index set $I_{k} \subseteq L$ and a number $\lambda_{k}>0$.

{\em Step 2:}  Take
$p^{k}=\pi_{Y_{(k)}}[y^{k-1}+\lambda_{k}(Ax^{k-1}-b)]$.

{\em Step 3:}  Take
$x^{k}={\rm argmin}\{f(x)+  \langle p^{k}, A x-b \rangle +0.5\lambda_{k}^{-1} \|x-x^{k-1}\|^{2} \ | \ x \in X\}$.

{\em Step 4:}  Take
$y^{k}=\pi_{Y_{(k)}}[y^{k-1}+\lambda_{k}(Ax^{k}-b)]$.
Set $k=k+1$ and go to Step 1.
\medskip

First we observe that
$$
p^{k}={\rm argmin}\{-\mathcal{L}(x^{k-1},p)+0.5\lambda_{k}^{-1} \|p-y^{k-1}\|^{2}  \ | \ p \in Y_{(k)}\}
$$
and
$$
y^{k}={\rm argmin}\{-\mathcal{L}(x^{k},y)+0.5\lambda_{k}^{-1} \|y-y^{k-1}\|^{2}  \ | \ y \in Y_{(k)}\}.
$$
Therefore, each iteration of (PDM) involves two projection (proximal) steps in the dual variable $y$
and one proximal step in the primal variable $x$. The point $w^{k}=(x^{k},y^{k})$
belongs to $X\times Y_{(k)}$. The next two properties follow the
usual substantiation schemes for this method; see \cite{Ant88} and also  \cite{CT94}.

\begin{lemma} \label{lm:3.1} Suppose $U $ is a closed convex set
in a finite-dimensional space $\mathbb{E}$, $\varphi : \mathbb{E} \to \mathbb{R}$ is a
convex function, $u$ is a point in $\mathbb{E}$. If
$$
\mu (z)=\varphi(z)+ 0.5\lambda ^{-1} \|z-u\|^{2}, \ \lambda>0,
$$
and
$$
v={\rm argmin}\{\mu (z) \ | \ z \in U\},
$$
then
\begin{equation} \label{eq:3.1}
 2\lambda  \{\varphi(v)-\varphi(z)\} \leq \|z-u\|^{2}-\|z-v\|^{2}-\|v-u\|^{2} \quad \forall z\in U.
\end{equation}
\end{lemma}
{\bf Proof.} Since the function $\mu$ is strongly convex
with constant $\lambda ^{-1}$, we have
$$
\mu (z)-\mu (v)\geq  0.5\lambda ^{-1} \|z-v\|^{2} \quad \forall z\in U.
$$
This inequality gives (\ref{eq:3.1}).
\QED


\begin{proposition} \label{pro:3.1}
Suppose that assumptions (A1)--(A2) are fulfilled.
For any pair $w^{*}=(x^{*},y^{*}) \in D^{*}_{(k)}\times Y^{*}_{(k)}$ we have
\begin{eqnarray}
\displaystyle
\|w^{k}-w^{*}\|^{2} && \leq  \|w^{k-1}-w^{*}\|^{2}-\|p^{k}-y^{k}\|^{2}-\|p^{k}-y^{k-1}\|^{2} -\|x^{k}-x^{k-1}\|^{2} \nonumber\\
  && +2 \lambda_{k} \langle y^{k}-p^{k},A(x^{k}-x^{k-1}) \rangle \nonumber\\
 && =  \|w^{k-1}-w^{*}\|^{2}-\|p^{k}-y^{k}\|^{2}-\|p^{k}-y^{k-1}\|^{2} -\|x^{k}-x^{k-1}\|^{2} \nonumber\\
  && +2 \lambda_{k}^{2} \|A_{(k)}(x^{k}-x^{k-1})\|^{2}. \label{eq:3.2}
\end{eqnarray}
\end{proposition}
{\bf Proof.} Choose any $w^{*}=(x^{*},y^{*}) \in D^{*}_{(k)}\times Y^{*}_{(k)}$.
Setting $\varphi(z)=\mathcal{L}(z,p^{k})$, $\lambda=\lambda_{k}$, $U=X$, $u=x^{k-1}$, $v=x^{k}$,
and $z=x^{*}$ in  (\ref{eq:3.1}) gives
$$
2\lambda_{k} \{\mathcal{L}(x^{k},p^{k})-\mathcal{L}(x^{*},p^{k})\} \leq \|x^{*}-x^{k-1}\|^{2}-\|x^{*}-x^{k}\|^{2}-\|x^{k}-x^{k-1}\|^{2}.
$$
Also, using  (\ref{eq:2.5}) with $I=I_{k}$, $x=x^{k}$,
and $y=p^{k}$ gives
$$
2\lambda_{k} \{\mathcal{L}(x^{*},p^{k})-\mathcal{L}(x^{k},y^{*})\} \leq 0.
$$
Adding these inequalities, we obtain
\begin{equation} \label{eq:3.3}
\|x^{k}-x^{*}\|^{2}  \leq  \|x^{k-1}-x^{*}\|^{2} -\|x^{k}-x^{k-1}\|^{2}
 +2 \lambda_{k} \langle p^{k}-y^{*},Ax^{k}-b \rangle .
\end{equation}
On the other hand, setting $\varphi(z)=-\mathcal{L}(x^{k-1},z)$, $\lambda=\lambda_{k}$, $U=Y_{(k)}$, $u=y^{k-1}$, $v=p^{k}$,
and $z=y^{k}$ in  (\ref{eq:3.1}) gives
$$
2\lambda_{k} \{\mathcal{L}(x^{k-1},y^{k})-\mathcal{L}(x^{k-1},p^{k})\} \leq \|y^{k}-y^{k-1}\|^{2}-\|p^{k}-y^{k}\|^{2}-\|p^{k}-y^{k-1}\|^{2}.
$$
Next,
setting $\varphi(z)=-\mathcal{L}(x^{k},z)$, $\lambda=\lambda_{k}$, $U=Y_{(k)}$, $u=y^{k-1}$, $v=y^{k}$,
and $z=y^{*}$ in  (\ref{eq:3.1}) gives
$$
2\lambda_{k} \{\mathcal{L}(x^{k-1},y^{*})-\mathcal{L}(x^{k},y^{k})\} \leq \|y^{*}-y^{k-1}\|^{2}-\|y^{*}-y^{k}\|^{2}-\|y^{k}-y^{k-1}\|^{2}.
$$
Adding these inequalities, we obtain
\begin{eqnarray}
\|y^{k}-y^{*}\|^{2}  && \leq  \|y^{k-1}-y^{*}\|^{2} -\|p^{k}-y^{k}\|^{2}-\|p^{k}-y^{k-1}\|^{2} \nonumber\\
 && -2 \lambda_{k} \{ \langle y^{*}-y^{k},Ax^{k}-b \rangle +\langle y^{k}-p^{k},Ax^{k-1}-b \rangle \}. \label{eq:3.4a}
\end{eqnarray}
Now adding (\ref{eq:3.3})  and (\ref{eq:3.4a}) gives the first inequality in (\ref{eq:3.2}).
Since
$$
\langle y^{k}-p^{k},A(x^{k}-x^{k-1}) \rangle =\lambda_{k} \|A_{(k)}(x^{k}-x^{k-1})\|^{2},
$$
we conclude also that the second relation in (\ref{eq:3.2}) holds true.
\QED

Now we can indicate conditions that provide basic convergence properties.


\begin{theorem} \label{thm:3.1}
Suppose that assumptions (A1)--(A2) are fulfilled,
\begin{equation} \label{eq:3.4}
\bigcap \limits^{\infty} _{k=1} W^{*}_{(k)} \neq \varnothing ,
\end{equation}
the sequence $\{\lambda_{k}\}$ satisfies the condition:
\begin{equation} \label{eq:3.5}
\lambda_k \in [\tau,\sqrt{(1-\tau)}/(\sqrt{2}\|A_{(k)}\|)]
\end{equation}
for some $\tau \in (0,1)$.
Then:
\begin{description}
  \item{(i)} the sequence $\{w^{k}\}$ has limit points,
  \item{(ii)} each of these limit points is a solution of problem (\ref{eq:2.5}) for some  $I \subseteq L$,
  \item{(iii)} for any limit point $\bar w$ of $\{w^{k}\}$ such that
$$
 \bar w \in \bigcap \limits^{\infty} _{k=1} W^{*}_{(k)}
$$
it holds that
\begin{equation} \label{eq:3.5a}
\lim \limits_{k\rightarrow \infty }w^{k}=\bar w.
\end{equation}
\end{description}
 \end{theorem}
{\bf Proof.} Take any point
$$
w^{*} \in \bigcap \limits^{\infty} _{k=1} W^{*}_{(k)}.
$$
Then from (\ref{eq:3.2}) and (\ref{eq:3.5}) we have
\begin{equation} \label{eq:3.6}
\|w^{k}-w^{*}\|^{2} \leq  \|w^{k-1}-w^{*}\|^{2}-\|p^{k}-y^{k}\|^{2}-\|p^{k}-y^{k-1}\|^{2} -\tau\|x^{k}-x^{k-1}\|^{2}
\end{equation}
for $k=1,2, \ldots$ Hence, the sequence $\{w^{k}\}$ is bounded and has
limit points, i.e. part (i) is true. Besides, (\ref{eq:3.6}) gives
\begin{equation} \label{eq:3.8}
\lim \limits_{k\rightarrow \infty }\|w^{k}- w^{*}\|=\sigma \geq 0
\end{equation}
and
\begin{equation} \label{eq:3.9}
\lim \limits_{k\rightarrow \infty }\|p^{k}-y^{k}\|=
\lim \limits_{k\rightarrow \infty }\|p^{k}-y^{k-1}\|=\lim \limits_{k\rightarrow \infty }\|x^{k}-x^{k-1}\|=0,
\end{equation}
hence
\begin{equation} \label{eq:3.10}
\lim \limits_{k\rightarrow \infty }\|y^{k}-y^{k-1}\|=0.
\end{equation}

Let $\bar w=(\bar x,\bar y)$ be an arbitrary limit point of $\{w^{k}\}$,
i.e.
$$
\bar w = \lim \limits_{{s} \rightarrow \infty } w^{k_{s}}.
$$
Then there exists  $J \subseteq L$ such that  $J = I_{k_{s}}$ for infinitely many times.
Without loss of generality we can suppose that $J = I_{k_{s}}$ for any $s$.
Then $w^{k_{s}}=(x^{k_{s}},y^{k_{s}}) \in X\times Y_{J}$ for any $s$, hence
$\bar w=(\bar x,\bar y)\in X\times Y_{J}$.
Setting $\varphi(z)=\mathcal{L}(z,p^{k})$, $\lambda=\lambda_{k}$, $U=X$, $u=x^{k-1}$, $v=x^{k}$,
and $z=x \in X$ in  (\ref{eq:3.1}) gives
$$
2\lambda_{k} \{\mathcal{L}(x^{k},p^{k})-\mathcal{L}(x,p^{k})\} \leq \|x-x^{k-1}\|^{2}-\|x-x^{k}\|^{2}-\|x^{k}-x^{k-1}\|^{2}.
$$
Taking the limit $k=k_{s} \rightarrow \infty $ due to (\ref{eq:3.9})--(\ref{eq:3.10}) gives
\begin{equation} \label{eq:3.11}
\mathcal{L}(\bar x,\bar y)-\mathcal{L}(x,\bar y) \leq 0.
\end{equation}
Also, setting $\varphi(z)=-\mathcal{L}(x^{k},z)$, $\lambda=\lambda_{k}$, $U=Y_{J}$, $u=y^{k-1}$, $v=y^{k}$,
and $z=y \in Y_{J}$ in  (\ref{eq:3.1}) gives
$$
2\lambda_{k} \{\mathcal{L}(x^{k},y)-\mathcal{L}(x^{k},p^{k})\} \leq \|y^{k-1}-y\|^{2}-\|y^{k}-y\|^{2}-\|y^{k}-y^{k-1}\|^{2}.
$$
Taking the limit $k=k_{s} \rightarrow \infty $ due to (\ref{eq:3.9})--(\ref{eq:3.10}) gives
\begin{equation} \label{eq:3.12}
\mathcal{L}(\bar x, y)-\mathcal{L}(\bar x,\bar y) \leq 0.
\end{equation}
It follows from (\ref{eq:3.11}) and (\ref{eq:3.12}) that $\bar w=(\bar x,\bar y)\in  W^{*}_{J}=D^{*}_{J}\times Y^{*}_{J}$.
Hence, part (ii) is also true.

Next, if
$$
 \bar w \in \bigcap \limits^{\infty} _{k=1} W^{*}_{(k)},
$$
 we can set $w^{*}=\bar w $ in (\ref{eq:3.8}). However, now $\sigma = 0$,
 which gives (\ref{eq:3.5a}) and part (iii) is true.  \QED

 These properties enable us to establish convergence
 to a solution under suitable conditions.


\begin{theorem} \label{thm:3.2}
Suppose that assumptions (A1)--(A2) are fulfilled,
the sequence $\{\lambda_{k}\}$ satisfies condition (\ref{eq:3.5}) for some $\tau \in (0,1)$.
\begin{description}
  \item{(i)} If there exists a nonempty basic index set $I \subseteq L$ such  that $I \subseteq I_{k}$,
  then  the sequence $\{w^{k}\}$ has limit points and each of these limit points
  belongs to  $ W^{*}$.
  \item{(ii)} If in addition $I = I_{k_{s}}$ for some infinite subsequence $\{w^{k_{s}}\}$, then
\begin{equation} \label{eq:3.13}
\lim \limits_{k\rightarrow \infty }w^{k}=w^{*} \in W^{*}.
\end{equation}
\end{description}
 \end{theorem}
{\bf Proof.} By definition, the sets $ W^{*}_{I}$ and $ W^{*}$ are now nonempty.
Due to Proposition \ref{pro:2.1}, $W^{*}_{I} \subseteq W^{*}_{(k)}$, hence
condition (\ref{eq:3.4}) holds. Then the sequence $\{w^{k}\}$ has limit points due to
Theorem \ref{thm:3.1} (i). Also,  there exists
$J \subseteq L$ such that  $I \subseteq J = I_{k_{s}}$ for infinitely many times.
But now $ J$ is a nonempty basic index set, hence  $W^{*}_{J} \subseteq W^{*}$.
Following the lines of part (ii) of Theorem \ref{thm:3.1}, we obtain that any
limit point of  $\{w^{k_{s}}\}$ will belong to  $ W^{*}_{J}\subseteq W^{*}$. Therefore,
part (i) is true. In case (ii) we have similarly that any
limit point $w^{*}$ of  $\{w^{k_{s}}\}$ will belong to  $ W^{*}_{I}$, but
$$
w^{*} \in \bigcap \limits^{\infty} _{k=1} W^{*}_{(k)}.
$$
The result now follows from Theorem \ref{thm:3.1} (iii). \QED


\begin{theorem} \label{thm:3.3}
Suppose that assumptions (A1)--(A2) are fulfilled, $F^{*} \bigcap  F_{L} \neq \varnothing$,
the sequence $\{\lambda_{k}\}$ satisfies condition (\ref{eq:3.5}) for some $\tau \in (0,1)$.
Then:
\begin{description}
  \item{(i)} the sequence $\{w^{k}\}$ has limit points,
  \item{(ii)} if each $I_{k}$ is a basic index set,
   all the limit points of $\{w^{k}\}$   belong to  $ W^{*}$,
  \item{(iii)} if there exists a nonempty basic index set $I \subseteq L$
  such that $I \subseteq I_{k}$ and $I = I_{k_{s}}$ for some infinite subsequence $\{w^{k_{s}}\}$,
  the sequence $\{w^{k}\}$ converges to a point of $W^{*}$.
\end{description}
 \end{theorem}
{\bf Proof.} Due to Proposition \ref{pro:2.2}, we now have
$F^{*} \bigcap  F_{I} = D^{*}_{I}$, $D^{*}_{I}\neq \varnothing$,
  and $\mathbf{0} \in  Y^{*}_{I}$ for any $I \subseteq L$. It follows that
$$
\left\{ F^{*} \bigcap  F_{L}\right\} \times \{\mathbf{0}\} \subseteq \bigcap \limits^{\infty} _{k=1} W^{*}_{(k)}.
$$
Therefore, (\ref{eq:3.4}) holds and assertion (i) follows from Theorem \ref{thm:3.1} (i).
Following the lines of part (ii) of Theorem \ref{thm:3.1}, we obtain that any
limit point of  $\{w^{k_{s}}\}$ will belong to  $ W^{*}_{J}\subseteq W^{*}$ where
$J$ is a nonempty basic index set. Therefore, assertion (ii) is also true.
Assertion (iii) clearly follows from Theorem \ref{thm:3.2}.\QED

The conditions of part (ii) of Theorem \ref{thm:3.2} are satisfied if for instance
we take the rule $I_{k} \subseteq I_{k+1}$ or $I_{k+1} \subseteq I_{k}$ for index sets.
These rules can be also applied in part (iii) of Theorem \ref{thm:3.3}.
In all the above theorems we utilized some conditions that must hold for each iterate $k$.
Obviously, all the assertions of the theorems will be true if we require for the same
conditions to hold only for $k \geq k'$ where  $k'$ is some fixed number.


\section{Primal-dual method for multi-agent optimization problems}\label{sc:4}

We now describe a specialization of the proposed approach
to the multi-agent optimization problem
\begin{equation} \label{eq:4.1}
 \min \to \left\{\sum \limits^{m} _{i=1} f_{i}(v) \ \vrule \ \bigcap \limits^{m} _{i=1} X_{i}  \right\},
\end{equation}
where $m$ is the number of agents (units) in the system. That is,
the information about the function $f_{i}$ and set $ X_{i}$ is known only to the $i$-th agent
and may be unknown even to its neighbours. Besides, it is usually supposed that
the agents are joined by some transmission links for information exchange
so that the system is usually a connected network, whose topology
may vary from time to time. This decentralized system has to find a concordant solution
defined by  (\ref{eq:4.1}).

For this reason, we replace
(\ref{eq:4.1}) with the family of optimization problems of the form
\begin{equation} \label{eq:4.2}
 \min \limits _{x \in D_{I}} \to f(x) = \sum \limits^{m} _{i=1} f_{i}(x_{i}),
\end{equation}
where $x=(x_{i})_{i=1, \ldots, m} \in \mathbb{R}^{mn}$, i.e. $x^{\top}=(x^{\top}_{1}, \dots,x^{\top}_{m})$,
$x_{i}= (x_{i1}, \dots,x_{in})^{\top}$  for $i=1, \dots, m$,
\begin{equation} \label{eq:4.3}
 D_{I} = X \bigcap F_{I}, \ X=X_{1}\times \dots \times X_{m}=\prod \limits_{i=1}^{m} X_{i}, \
 X_{i}\subseteq \mathbb{R}^{n}, \ i=1, \dots, m;
\end{equation}
the set $F_{I}$ describes the information exchange scheme within the current topology of the communication network,
and $I$ is the index set of arcs of the corresponding oriented graph.
More precisely, the maximal (full) communication network with non-oriented edges denoted by $\mathcal{F}$ corresponds to
the set
$$
\tilde F=\left\{x \in \mathbb{R}^{mn} \ | \ x_{s}=x_{t}, \ s,t=1, \dots, m, \ s \neq t  \right\},
$$
i.e. each edge is  associated with two directions or equations ($x_{s}=x_{t}$ and $x_{t}=x_{s}$).
However, this definition of topology is superfluous. It seems more suitable to introduce
some other graph topology for writing the multi-agent optimization problem
in addition to the graph $\mathcal{F}$.
For this reason, we associate each pair of vertices (agents) $(s,t)$
to one oriented arc $i$, so that $L= \{1, \dots, l\}$ is the index set of all these arcs,
hence $l=m(m-1)/2$. That is, each arc $(s,t)$ is in fact used in both the directions in
the communication network $\mathcal{F}$, but we fix only one direction for definition of
the multi-agent optimization problem and obtain the graph $\mathcal{G}$.
Taking subsets $I \subseteq L$, we obtain various constraint sets
\begin{equation} \label{eq:4.4}
 F_{I} = \left\{x \in \mathbb{R}^{mn} \ | \ x_{s}-x_{t}=\mathbf{0}, \ i=(s,t)  \in I  \right\},
\end{equation}
corresponding to the oriented graphs $\mathcal{G}_{I}$ in the
the multi-agent optimization problem formulation. Replacing
the arcs in $\mathcal{G}_{I}$ with non-oriented edges, we obtain
the corresponding communication network $\mathcal{F}_{I}$ of the system.
It follows that
$\mathcal{F}=\mathcal{F}_{L}$, $\mathcal{G}=\mathcal{G}_{L}$,
and $F =F_{L}$. Next, for each arc $i=(s,t)$ we can define the $n \times mn$ sub-matrix
$$
A_{i}=\left( A_{i1} \cdots A_{im} \right),
$$
where
$$
A_{ij}= \left\{
\begin{array}{rl}
E, & \mbox{if} \ j=s, \\
-E, & \mbox{if} \ j=t, \\
\Theta, & \mbox{otherwise},
\end{array}
 \right.
$$
$E$ is the $n \times n$ unit matrix, $\Theta$ is the $n \times n$ zero matrix.
Then clearly
$$
F_{I} = \left\{x \in \mathbb{R}^{mn} \ | \ A_{I} x=\mathbf{0}  \right\},
$$
where
$$
A_{I}=\left( \{A^{\top}_{i}\} _{i \in I}\right)^{\top},
$$
which corresponds to the definition in (\ref{eq:2.3})
for $b_{I}=\mathbf{0}$ and any $I \subseteq L$, hence we can set $A=A_{L}$.
Therefore, our problem  (\ref{eq:4.2})--(\ref{eq:4.4})
corresponds to  (\ref{eq:2.3})--(\ref{eq:2.4}).

In what follows, we will use the following basic assumptions.
\begin{enumerate}
\item[{\bf (B1)}] For  each $i=1, \dots, m$, $X_{i}$ is a  convex and closed set in
$ \mathbb{R}^{n}$, $f_{i} : \mathbb{R}^{n} \to \mathbb{R}$ is a convex
function.
\item[{\bf (B2)}]
The set $ D^{*}=D^{*}_{L}$ is nonempty.
\end{enumerate}

These assumptions imply (A1)--(A2). If the graph $\mathcal{F}_{I}$
for some $I \subseteq L$ is connected,  then $I$ a basic index set.
Now we present an implementation of Method (PDM) for
the multi-agent optimization problem  (\ref{eq:4.2})--(\ref{eq:4.4}),
where each agent (or unit) receives information only from
its neighbours. Given an oriented graph $\mathcal{G}_{I}$
and an agent $j$, we denote
by $\mathcal{N}^{+}_{I}(j)$ and $\mathcal{N}^{-}_{I}(j)$ the sets of
incoming and outgoing arcs at $j$. Since many
oriented graphs $\mathcal{G}_{I}$ are associated with the same
graph $\mathcal{F}_{I}$, we suppose that agent $j$ is responsible
for calculation of the current values of the primal variable
$x_{j}$ and all the dual variables $y_{i}$ and $p_{i}$ such that
$i \in \mathcal{N}^{-}_{I}(j)$. That is, we will fix the
oriented graph $\mathcal{G}$ and its subgraphs $\mathcal{G}_{I}$ such that
agent $j$ is associated with all the outgoing arcs for vertex $j$.
The general Lagrange function for problems  (\ref{eq:4.2})--(\ref{eq:4.4}) is written as follows:
\begin{eqnarray}
 \mathcal{L}(x,y) &= & f(x)+  \langle y, A x\rangle
 =  \sum \limits _{j\in M} f_{j}(x_{j})
           +  \sum \limits _{i\in L} \langle y_{i}, A_{i} x\rangle  \label{eq:4.5}\\
 &= & \sum \limits _{j\in M} \left\{ f_{j}(x_{j})
 +  \sum \limits _{i\in \mathcal{N}^{-}_{L}(j)} \langle y_{i}, x_{j}\rangle -
 \sum \limits _{i\in \mathcal{N}^{+}_{L}(j)} \langle y_{i}, x_{j}\rangle \right\}. \nonumber
\end{eqnarray}
The saddle point problems are defined in (\ref{eq:2.5}).
As in Section \ref{sc:3}, for simplicity we will write $Y_{(k)}=Y_{I_{k}}$,
$Y^{*}_{(k)}=Y^{*}_{I_{k}}$, etc.

\medskip
\noindent {\bf Method (PDMI).}
At the beginning, the agents choose the communication topology by
choosing the active arc index set $I_{0} \subseteq L$. Next,
each $s$-th agent chooses $x^{0}_{s}$ and $y^{0}_{i}$ for
$i \in \mathcal{N}^{-}_{(0)}(s)$ and reports these values
to its neighbours. This means that $y^{0}_{i}=\mathbf{0}$ for $i \notin I_{0}$.

 At the $k$-th iteration, $k=1,2,\ldots$,
 each $s$-th agent has the values $x^{k-1}_{s}$ and $y^{k-1}_{i}$,
$i \in \mathcal{N}^{-}_{(k-1)}(s)$, and the same values of its neighbours.
The agents choose the current communication topology by
choosing the active arc index set $I_{k} \subseteq L$ and determine
the stepsize $\lambda_{k}$. This means that they set $y^{k}_{i}=\mathbf{0}$ for $i \notin I_{k}$.

{\em Step 1:}  Each $s$-th agent sets
\begin{equation} \label{eq:4.6}
p^{k}_{i}=y^{k-1}_{i}+\lambda_{k} (x_{s}^{k-1}-x_{t}^{k-1}) \quad  \forall i=(s,t), \ i \in \mathcal{N}^{-}_{(k)}(s).
\end{equation}
Then each $s$-th agent reports these values to its neighbours.

{\em Step 2:}  Each $s$-th agent calculates
$$
v_{s}^{k}=\sum \limits _{i\in \mathcal{N}^{-}_{(k)}(s)} p^{k}_{i} - \sum \limits _{i\in \mathcal{N}^{+}_{(k)}(s)}  p^{k}_{i}
$$
and
\begin{equation} \label{eq:4.7}
\displaystyle
x_{s}^{k}=\arg\min_{x_{s} \in X_{s}} \left\{f_{s}(x_{s})+ \langle v_{s}^{k}, x_{s}\rangle
     +0.5\lambda_{k}^{-1} \|x_{s}-x_{s}^{k-1}\|^{2} \right\}
\end{equation}
and reports this value to its neighbours.

{\em Step 3:}  Each $s$-th agent sets
\begin{equation} \label{eq:4.8}
y^{k}_{i}=y^{k-1}_{i}+\lambda_{k} (x_{s}^{k}-x_{t}^{k}) \quad  \forall i=(s,t), \ i \in \mathcal{N}^{-}_{(k)}(s).
\end{equation}
Then each $s$-th agent reports these values to its neighbours.
The  $k$-th iteration is complete.
\medskip

We observe that the agents do not store the dual variables
related to the inactive arcs, i.e. $y^{k}_{i}=\mathbf{0}$ for $i \notin I_{k}$.
If some arc $i =(s,t) \notin I_{k-1}$ becomes active at the $k$-th iteration,
i.e. $i \in I_{k}$, then agent $s$ simply sets $y^{k-1}_{i}=\mathbf{0}$.

Due to (\ref{eq:4.5}), relations (\ref{eq:4.6})--(\ref{eq:4.8}) correspond to
Steps 2--4 of (PDM), respectively. Hence,  the convergence properties
of (PDMI) will follow directly from
Theorems \ref{thm:3.2} and \ref{thm:3.3}.


\begin{corollary} \label{cor:4.1}
Suppose that assumptions (B1)--(B2) are fulfilled,
the sequence $\{\lambda_{k}\}$ satisfies condition (\ref{eq:3.5}) for some $\tau \in (0,1)$.
\begin{description}
  \item{(i)} If there exists a nonempty basic index set $I \subseteq L$ such  that $I \subseteq I_{k}$,
  then  the sequence $\{w^{k}\}$, $w^{k}=(x^{k},y^{k})$, generated by (PDMI) has limit points and each of these limit points
  belongs to  $ W^{*}$.
  \item{(ii)} If in addition $I = I_{k_{s}}$ for some infinite subsequence $\{w^{k_{s}}\}$, then
(\ref{eq:3.13}) holds.
\end{description}
\end{corollary}


\begin{corollary} \label{cor:4.2}
Suppose that assumptions (B1)--(B2) are fulfilled,
the sequence $\{\lambda_{k}\}$ satisfies condition (\ref{eq:3.5}) for some $\tau \in (0,1)$, $F^{*} \bigcap  F_{L} \neq \varnothing$.
Then:
\begin{description}
  \item{(i)} the sequence $\{w^{k}\}$, $w^{k}=(x^{k},y^{k})$, generated by (PDMI) has limit points,
  \item{(ii)} if each $I_{k}$ is a basic index set,
   all the limit points of $\{w^{k}\}$   belong to  $ W^{*}$,
  \item{(iii)} if there exists a nonempty basic index set $I \subseteq L$
  such that $I \subseteq I_{k}$ and $I = I_{k_{s}}$ for some infinite subsequence $\{w^{k_{s}}\}$,
  the sequence $\{w^{k}\}$ converges to a point of $W^{*}$.
\end{description}
\end{corollary}

Convergence of (PDMI) requires for all the agents to choose the stepsize $\lambda_{k}$
in accordance with (\ref{eq:3.5}), hence they have to evaluate the norm $\|A_{(k)}\|$
at the $k$-th iteration. Fix some $I \subseteq L$,  then
$$
A_{I}^{\top} A_{I}=  H_{I} \otimes E,
$$
where $H_{I}$ is the Kirchhoff matrix of the graph $\mathcal{F}_{I}$,
$\otimes$ denotes the Kronecker product of matrices. Application of
the Gershgorin theorem (see Theorem 5 in \cite[Chapter XIV]{Gan66}) gives
$$
\|A_{I}\|=\sqrt{\|H_{I}\|} \leq \sqrt{2 d(\mathcal{F}_{I})},
$$
where $d(\mathcal{F}_{I})$ is the maximal vertex degree of the graph $\mathcal{F}_{I}$.
There exist more precise estimates for some special classes of graphs; see e.g.
\cite{LZ97,Pan02}. Together with (\ref{eq:3.5}) we obtain the bound
\begin{equation} \label{eq:4.9}
\lambda_k \in \left[\tau,0.5 \sqrt{(1-\tau)/d(\mathcal{F}_{(k)})}\right]
\end{equation}
for some $\tau \in (0,1)$. In case of varying topology the separate agents
may meet difficulties in evaluation of $d(\mathcal{F}_{(k)})$ since the graph
then may be non-regular. The concordant value of $\lambda=\lambda_k$
satisfying (\ref{eq:4.9}) can be obtained by determining some upper bound for
$d(\mathcal{F}_{(k)})$. It seems suitable to apply  the following strategy.
First we choose the fixed topology that corresponds to an arc index set $J \subset L$
so that it gives the connected graph $\mathcal{F}_{J}$ and $J \subseteq I_{k}$ for any $k$.
This means that all the arcs in $J$ remain always active. The status of
the other arcs may vary, but the maximal vertex degree of the graph $\mathcal{F}_{(k)}$
can not exceed some fixed number $v$. Then each agent can take
$\lambda= 0.5 \sqrt{(1-\tau)/v}$ and the assumptions of
Corollary \ref{cor:4.1} (i) and Corollary \ref{cor:4.2} (i)--(ii) on the choice of parameters hold.

We now give a natural example of problem  (\ref{eq:4.2})--(\ref{eq:4.4})
such that $F^{*} \bigcap  F_{L} \neq \varnothing$. Namely, set
$ X_{i}=\mathbb{R}^{n}$, $f_{i}(v)=(1/p)(\max\{h_{i}(v),0\})^{p}$, $p \geq 1$
for $i=1, \ldots, m$. Then (\ref{eq:4.2})--(\ref{eq:4.4}) corresponds to a penalized
problem for finding a point of the set
$$
\tilde V=\left\{u \in \mathbb{R}^{n} \ | \ h_{i}(u) \leq 0, \ i=1, \dots, m \right\}.
$$
If $\tilde V \neq \varnothing$, then clearly $F^{*} \bigcap  F_{L} \neq \varnothing$,
which gives stronger convergence properties.

It should be noticed that primal-dual methods are usually applied to
large-scale convex optimization problems with binding constraints in order to
keep the decomposability properties. However, the streamlined
primal-dual gradient projection method requires strengthened
assumptions. Utilization of extrapolation steps enables one to attain
convergence under custom convex-concavity; see \cite{AS58}.
These methods admit a fixed positive stepsize that yields a linear rate of
convergence; see e.g. \cite[Chapter VI]{GT89} and the references therein.
However, replacing projections with proximal steps also will
enhance convergence, besides the method becomes applicable to non-smooth problems.
This primal-dual method with proximal steps was proposed in \cite{Ant88}.
Similar methods were described in \cite{CT94,EZC10}. It should be also noticed that
known iterative methods for multi-agent optimization problems with changing communication
topology are based on different conditions; see e.g. \cite{NO15,AH16}.


\end{document}